\documentclass[10pt]{amsart}
\usepackage[cp1251]{inputenc}
\usepackage[english]{babel}
\usepackage{amsmath}
\usepackage{amssymb}
\usepackage{amsfonts}
\usepackage{srcltx} % Ї ЄҐв ¤«п ®Ўа в­®Ј® Ї®ЁбЄ  (Ё§ DVI ў вҐЄбв®ўл©)
\usepackage[dvips]{graphicx}
\newtheorem{lemma}{Lemma}
\newtheorem{theorem}{Theorem}

\newcommand{\M}{{\mathcal{M}}}

\newcommand{\F}{{\mathbb{F}}}

\newcommand{\zero}{{\mathbf{0}}}

\newcommand{\iso}{\mbox{Iso}}
\newcommand{\sym}{\mbox{Sym}}

\begin{document}

\title{Ranks of propelinear perfect binary codes}
\author{{G. K. Guskov, I. Yu. Mogilnykh,  F. I. Solov'eva}}%
\address{George K. Guskov
\newline\hphantom{iii} Sobolev Institute of Mathematics,
\newline\hphantom{iii} pr. ac. Koptyuga 4,
%%%# \newline\hphantom{iii} Novosibirsk State University,%
%%%# \newline\hphantom{iii} Pirogova street~2,
\newline\hphantom{iii} 630090, Novosibirsk, Russia}%
\email{m1lesnsk@gmail.com}%

\address{Ivan Yu. Mogilnykh
\newline\hphantom{iii} Sobolev Institute of Mathematics,
\newline\hphantom{iii} pr. ac. Koptyuga 4,
\newline\hphantom{iii} Novosibirsk State University,%
\newline\hphantom{iii} Pirogova street~2,
\newline\hphantom{iii} 630090, Novosibirsk, Russia}%
\email{ivmog84@gmail.com}%

\address{Faina I. Solov'eva
\newline\hphantom{iii} Sobolev Institute of Mathematics,
\newline\hphantom{iii} pr. ac. Koptyuga 4,
\newline\hphantom{iii} Novosibirsk State University,%
\newline\hphantom{iii} Pirogova street~2,
\newline\hphantom{iii} 630090, Novosibirsk, Russia}%
\email{sol@math.nsc.ru}%

%\thanks{\sc Ivanov, I.I., On some problems of commutative algebra}
\thanks{\copyright \ 2012 G. K. Guskov, I. Yu. Mogilnykh,  F. I. Solov'eva}
\thanks{\rm The second author was supported by the Grant of
the President of the Russian Federation for Young Russian
Researchers (project no. MK-1700.2011.1) and by the Grants RFBR
%the Russian Foundation for Basic Research (project nos.
 12-01-00448- , 12-01-31098 and 10-01-00616-a. The
work of the third author was partially supported by Grants RFBR
10-01-00424-a and 12-01-00631-a.}
%\thanks{\it Џ®бвгЇЁ«  1 ­®пЎап 2012 Ј., ®ЇгЎ«ЁЄ®ў ­  1 ¤ҐЄ Ўап 2012 Ј.}%

%\top

%\vspace{1cm}

\maketitle

\begin{quote}
{\small \noindent{\sc Abstract. } It is proven that for any
numbers $n=2^m-1, m\geq 4$ and $r$, such that $n-\log(n+1)\leq r
\leq n$ excluding $n=r=63$, $n=127$, $r\in\{126,127\}$ and
$n=r=2047$ there exists a propelinear perfect binary code of
length $n$ and rank $r$.
\medskip

\noindent{\bf Keywords:} propelinear perfect binary codes, rank,
transitive codes}
\end{quote}

\section{Introduction}

Denote by $\F^n$ a vector space of dimension $n$ over the Galois
field $GF(2)$ with respect to  the Hamming distance. The  {\em
Hamming distance} $d(u,v)$ between two vectors $u,v\in\F^n$ is
defined as the number of coordinates in which $u$ and $v$ differ.
Any isometry of $\F^n$ is given by a coordinate permutation and a
translation.  We denote by $\iso(\F^n)$ the group of all
isometries of $\F^n$:
$$
\iso(\F^n)=\{(v,\pi)\mid v\in\F^n,\pi\in  S_n\},
$$
where $ S_n$ denotes the symmetric group of degree $n$ and
$
(v,\pi)(x)=v+\pi(x)$ for any $x\in \F^n.
$
The group operation in $\iso(\F^n)$ is the composition
$
(u,\pi)\circ (v,\tau)=(u+\pi(v),\pi\tau)$ for all $(u,\pi),(v,\tau)\in \iso(\F^n).
$
Here, and throughout the entire paper, we use
$\pi\tau(x)=\pi(\tau(x))$ for $x\in\F^n$.

An arbitrary subset of $\F^n$ is called a binary code of length
$n$. The {\it minimum distance} of a code $C$ is the minimum value
of the Hamming distance between any two different codewords from
$C$. Two codes $C$ and $D$ are said to be {\em equivalent} if
$C=\phi(D)$, for some isometry $\phi$ of $\F^n$. By $\sym(C)$ we
denote the group of all coordinate permutations that fix the code
$C$ set-wise and call it the {\it symmetry group} of $C$. By
$\iso(C)$ we denote the group of all isometries of $\F^n$ fixing
the code $C$ set-wise, and we call it the {\it automorphism group}
of $C$. Note that in some papers, code automorphisms are defined
as coordinate permutations fixing the code set-wise.

A code $C$ is called  {\it single-error-correcting
perfect} (or  perfect, for the sake of brevity)
 if for any vector $x\in \F^n$ there exists exactly one vector $y
\in C$ such that $d(x,y)\leq 1$. It is well known that such codes
exist  if and only if $n=2^m-1, m\geq 1$. For any
$n=2^m-1, m\geq 1$, there is exactly one, up to equivalence,
linear perfect code of length $n$ and it is called the {\it
Hamming code}.

 Throughout  the paper we
assume that $C\in \F^n$ is a perfect code of length $n$ containing
the all-zero vector  $\zero^n$ with  $n$ coordinates. For such a
code $C$,  its kernel $K$ is defined as the set of all
codewords that leave $C$ invariant under translation, that is,
$$
K=\{x\in C \mid x+C=C\}.
$$
The kernel $K$ of $C$  is a linear subspace of $\F^n$ and the code
$C$ is a union of cosets of $K$. Rank $rank(C)$ of a code $C$ is
the dimension of 
the linear span $<C>$. By $e_i$ we denote a
vector of weight $1$ having unit in $i$th coordinate position.

A code $C$ is called {\em transitive} if $\iso(C)$ acts
transitively on $C$.

Let $\Pi$ be a mapping of the codewords from
 $C$ into the admissible
permutations: $x\mapsto \pi_{x}$: $(x,\pi_x) \in \iso(C)$, such
that $\pi_{(x,\pi_{x})y}=\pi_{x}\pi_{y}$. Then we can define a
group operation on $C$: $$x\star y=(x,\pi_x)y.$$ A code equipped
with the operation defined above is called {\it a propelinear
structure} on $C$ and is denoted by $(C,\Pi,\star)$ (simply
$(C,\star)$ if we do not need any information on $\Pi$). A code is
called {\it propelinear} if it has a propelinear structure.

It is easy to see that any propelinear code is
transitive. Transitive codes were constructed and studied
 in
\cite{S2004,S2005}. Propelinear codes were introduced in
\cite{Rif1} and investigated further in
\cite{Rif2,BMRS2012,BMRS2013}.  It is proven
 that perfect propelinear codes can be obtained by using the well known  Vasil'ev
construction, see  \cite{RPB}, and by the Mollard construction, see the proof in \cite{BMRS2012}.
 In \cite{BMRS2013} an
  exponential number of nonequivalent propelinear perfect codes having small ranks is presented.

 In this paper we solve the rank problem for propelinear perfect
 codes: all possible ranks of perfect codes are attainable by propelinear perfect codes,
 except full ranks for lengths 63, 127, 2047 and the rank 126 for
 codes of length 127.

\section{Propelinear full rank perfect codes  of lengths $15$ and $31$}

 Let us recall the Vasil'ev construction
\cite{V62}.  Let  $C$ be a
%!!
perfect binary code of length
$(n-1)/2$. Let $\lambda $ be any map from $C$ into the set
$\{0,1\}$  and $|x|=x_1+\dots +x_n$, where $x=(x_1,\ldots ,x_n)$,
$x_i\in \{0,1\} $. The code
\begin{equation} \label{Vas62} C^{n} =
\{(x+y,|x|+\lambda (y),x) \mid x \in \F^{(n-1)/2}, y \in
C\}\end{equation}   is called  {\it Vasil'ev perfect code}. Let
$(C,\star)$ be a propelinear structure on $C$, then a homomorphism
$\lambda$ from $(C,\star)$ into $Z_2$ is called a {\it propelinear
homomorphism} (or {\it propelinear function}).

\begin{theorem} (See \cite{RPB}) \label{VasilevPropelinear}
Let $(C,\star)$ be a
%!!
 propelinear structure on a perfect binary
code $C$ of length $n$, let $\lambda$ be a propelinear function
from the code $C$ into $Z_2$.
 Then the Vasil'ev code $C^{n}$ is propelinear perfect.
 \end{theorem}

Generally speaking, the problem of checking propelinearity of a
given transitive code is computationally hard. In \cite{BMRS2012}
we limited ourselves with normalized propelinear codes. Recall,
see \cite{BMRS2012}, that a propelinear structure $(C,\Pi,\star)$
is called a {\em normalized propelinear}
 if the permutations assigned to the codewords of the same coset of the kernel, coincide. Computer research
is carried out in a way that the number of possible candidates for
propelinear structures increases exponentially as the size of
kernels decreases by unity, meaning that codes of full rank seem
to be out of a computational reach
%% ўбҐ-в ЄЁ reach -ў­Ґ ¤®бвЁ¦Ё¬®бвЁ Є®¬ЇмовҐа 
 (as
 they have relatively small kernels). In order
to avoid this problem, we require codes to have trivial symmetry
groups. In this case, there is just one opportunity for a
assignment of permutations, in other words, $Aut(C)$ is acting
regularly on codewords of $C$,
 \cite{Rif3}. So, $C$ is a normalized propelinear
code and the following statement holds. %% ‡¤Ґбм Ўл«® Ї«®е®Ґ ЇаҐ¤«®¦Ґ­ЁҐ, Ї®¬Ґ­п«.

\begin{lemma}\label{trivsym}
A transitive code with trivial symmetry group is normalized
propelinear.
\end{lemma}

Among perfect codes of length 15 from the database \cite{OstPot}, we
found $44$ transitive codes with trivial symmetry groups,
 $39$ of
them having full rank and $5$ having rank $14$. Note that the
existence of propelinear perfect codes of length 15 of all
possible ranks, with the exception of full rank code,
%% we study the existence of codes, but not ranks
 was previously
shown in \cite{BMRS2012}.

\begin{lemma}\label{full rank for 15}
There is a propelinear perfect code of length 15 of any admissible
rank.
\end{lemma}

We give two more lemmas concerning Vasil'ev codes. Note that the
assigned permutations $\Pi(C)$ of the propelinear code $C$ of
length $(n-1)/2$ form a subgroup of $S_{(n-1)/2}$, see
\cite{BMRS2012}. Some of the homomorphisms of $C$ into $Z_2$ can
be described in terms of those of the group $\Pi(C)$.

\begin{lemma}\label{extension}
Let $(C,\Pi,\star)$ be a propelinear code. Any group homomorphism
 $\lambda'$ of $(\Pi(C),\circ)$ into $Z_2$ can be extended to
a propelinear homomorphism $\lambda$ of $(C,\Pi,\star)$ into $Z_2$
in the following way: $\lambda(x):=\lambda'(\pi_x)$.
\end{lemma}
{\bf Proof}. The structure-preserving property follows immediately
from the definition of a propelinear code:
$$\lambda(x\star y)=\lambda'(\pi_{x\star y})=\lambda'(\pi_{(x,\pi_{x})y})=\lambda'(\pi_x\pi_y)=\lambda'(\pi_x)+\lambda'(\pi_y)=\lambda(x)+\lambda(y).$$

\begin{lemma}\label{reduction}
Let $C^{n}$ be a
 code given by
 the  Vasil'ev construction
(\ref{Vas62}) with
 function $\lambda$. Then
 $rank(C^{n})=rank(\{(y,\lambda(y)): y \in C\})+(n-1)/2$ \,\, and

$rank(C)+(n-1)/2 \leq rank(C^n)\leq rank(C)+(n+1)/2$. \end{lemma}
{\bf Proof}. The basis of the linear span of $C^n$ can be chosen
in such a way that it contains vectors: $(x^{i},|x^{i}|,x^{i})$,
for vectors $\{x^{i}:i \in \{1, \ldots , (n-1)/2\}\}$ being a
basis of $F^{(n-1)/2}$. Obviously, the rank of
$\{(y,\lambda(y),{\bf 0}^{(n-1)/2}): y \in C\}$ is equal to that
of $\{(y,\lambda(y)): y \in C\}$.

 Depending on the function $\lambda$ the rank of the code $C^n$ is equal to
$rank(C)+(n+1)/2$ if the vector $e_{n+1}$ belongs to its span,
otherwise it is equal to $rank(C)+(n-1)/2$.

\begin{theorem} \label{full rank for 31}
There exists a  full rank normalized %% added normalized here.
 propelinear perfect binary
code of length $31$.
\end{theorem}
{\bf Proof}. Lemma \ref{full rank for 15} implies the existence of
propelinear perfect codes of length $15$ of full rank. In order to
construct a perfect code of length $31$ of full rank, another
computer search was carried out. As mentioned before, there are
exactly  39
 propelinear
 full rank  perfect codes of length 15 with trivial symmetry group. For each
of the codes we considered propelinear homomorphisms of special
type, i.e., satisfying Lemma \ref{extension} and looked at the
sizes of the ranks of the Vasil'ev codes of length 31 using Lemma
\ref{reduction}. Only three  of 39 codes (the numbers of these codes are 5584, 5844, 5823 from the database \cite{OstPot}) produce full rank
Vasil'ev codes of length 31. An interesting fact is that the
symmetry groups of the Steiner triple systems of the obtained
codes of length 31 are trivial,
 so the codes inherit the trivial symmetry group property.

\section{Rank problem}

In this section we solve the rank problem for propelinear perfect
codes using the results of the previous section as well as the Vasil'ev and
the Mollard constructions.
 Recall   the Mollard
construction  for binary codes. Let $C^t$ and $C^m$ be any two
perfect codes of lengths $t$ and $m$, respectively, containing
all-zero vectors.

Let $x=(x_{11},x_{12}, \ldots,x_{1m},
x_{21},\ldots,x_{2m},\ldots,x_{t1},\ldots,x_{tm}) \in\,\F^{tm}.$
 The
generalized parity-check functions $p_{1}(x)$ and $p_{2}(x)$ are
defined as $p_{1}(x)=(\sigma_{1},\sigma_{2},\ldots,\sigma_{t}) \in
\F^{t},$ $p_{2}(x)=(\sigma'_{1},\sigma'_{2},\ldots,\sigma'_{m})
\in \F^{m},$ where $\sigma_{i}=\sum_{j=1}^{m}x_{ij}$ and
$\sigma'_{j}=\sum_{i=1}^{t}x_{ij}$. Let $f$ be any function from
$C^t$ to $\F^m.$  The set
$$\M(C^{t},C^{m})=\{(x,y\, +
\,p_{1}(x),z\, + \,p_{2}(x)\, + \, f(y)) \mid
x\in\,\F^{tm},y\in\,C^t, z \in\,C^{m}\}$$ is a perfect binary
Mollard code of length $n=tm+t+m$, see \cite{Mollard}. Here the
abbreviation $\M(C^{t},C^m)$ indicates the lengths of initial
codes $C^t$ and $C^m$. It is clear that the codes with other
lengths $t'$ and $m'$ can also yield a perfect code
$\M(C^{t'},C^{m'})$ with the same parameters as the code
$\M(C^{t},C^{m})$,  both these codes  could coincide or  be
different, moreover, they could be nonequivalent.

\begin{theorem} (See \cite{BMRS2012}) \label{Mollard Propelinear} Let $C^t$ and $C^m$ be arbitrary
 propelinear perfect binary codes
of lengths $t$ and $m$, respectively. Let $f$ be a propelinear
homomorphism from $C^t$ to $\F^m.$ Then the Mollard code
$\M(C^{t},C^{m})$   is a propelinear perfect binary
 code of length $n=tm+t+m$, see \cite{BMRS2012}.
\end{theorem}

Further we consider the Mollard codes with the function $f\equiv {\bf 0}^m.$

\begin{lemma}(See \cite{S2005}) \label{ranks}
The perfect binary Mollard code $\M(C^{t},C^{m})$ of length
$n=tm+t+m$ with $f\equiv {\bf 0}^m$ has rank $tm + r(C^t) + r(C^m)$.
\end{lemma}

\begin{theorem}
For any $n=2^m-1, m\geq 4$ and  arbitrary $r$, satisfying
$n-\log(n+1)\leq r \leq n$ excluding cases of
$n=r=63$; $n=127$, $r\in\{126,127\}$ and $n=r=2047$,  there exists
a propelinear
 perfect
binary code of length $n$ and rank $r$.
\end{theorem}

{\bf Proof}. The proof is provided by  applying the Vasil'ev
construction for small $n$ and by induction applying the  Mollard
construction beginning with $n=2^8-1$. In order to make the
induction step working we need several initial steps.

By Lemma \ref{full rank for 15}  for $n=15$ we have propelinear
perfect codes of length 15 of all possible ranks.
% Table 1 in the paper \cite{BMRS2012} shows
%that there are propelinear perfect codes of length 15 with ranks
%$11$ (the Hamming code), $12$, $13$ and $14$;  the
%propelinear perfect code of length 15 with full rank $15$ is
%constructed in the previous section.

Using these propelinear codes of length $15$, Theorem
\ref{VasilevPropelinear} and Lemma \ref{reduction} setting the
function $\lambda\equiv 0$ we obtain propelinear perfect codes of
length $31$ having all possible ranks with the exception of full
rank. A full rank code we have by Theorem \ref{full rank for 31}.

Applying further the Vasil'ev construction with the function
$\lambda\equiv 0$ we obtain for $n=63$ propelinear perfect codes
of all possible ranks, except the full rank. For $n=127$ we start
with the obtained Vasil'ev perfect codes of length $63$ and again
by  the Vasil'ev construction with $\lambda\equiv 0$ we obtain
propelinear codes of length $127$ for all possible ranks with the
exceptions of codes of full rank and rank $126$.

Let us consider  the  Mollard codes

\begin{equation}\label{Mollard small}\M(C^{2^{4}-1},C^{2^4-1}), \,
\M(C^{2^{4}-1},C^{2^5-1}), \, \M(C^{2^{5}-1},C^{2^5-1})
\end{equation}
of lengths 255, 511 and 1023 respectively. From Lemma \ref{ranks}
varying the propelinear codes of different ranks of lengths 15 and
31, we get the propelinear Mollard codes (\ref{Mollard small}) for
each possible rank.

 In order to fulfill the
case $n=r=2^{11}-1=2047$ we have to construct the Mollard code
$\M(C^{2^4-1},C^{2^7-1})$ or $\M(C^{2^5-1},C^{2^6-1})$ from full
rank propelinear codes of length 63 or 127, which we do not have
(or we have to use another approach to construct such codes). But
as we see below the open cases do not influent on the process of
obtaining propelinear perfect codes of all possible ranks and all
admissible lengths $n\geq 2^{12}-1$.

Let the theorem be true and there exist propelinear perfect codes
of any rank for every  length
$$2^{4s}-1, \,\, 2^{4s+1}-1, \,\, 2^{4s+2}-1$$
 for $s\geq 2$.

Applying the Mollard construction to  these propelinear codes and
propelinear perfect codes of length $15$ or $31$ of different
ranks by Theorem \ref{Mollard Propelinear}  we obtain the
following four  perfect codes
\begin{equation}\label{Mollardformulas}\M(C^{2^{4s}-1},C^{2^4-1}), \,
\M(C^{2^{4s}-1},C^{2^5-1}), \,
\end{equation}
$$ \M(C^{2^{4s+1}-1},C^{2^5-1}), \,
\M(C^{2^{4s+2}-1},C^{2^5-1}),$$
of lengths
\begin{equation}\label{Mollard lengths}2^{4(s+1)}-1, \,\, 2^{4s+5}-1, \,\, 2^{4s+6}-1, \,\,
2^{4s+7}-1,\end{equation}  respectively.  From Lemma \ref{ranks}
we see that varying the codes of different ranks in the induction
hypotheses,  we obtain the Mollard codes (\ref{Mollardformulas})
for every length (\ref{Mollard lengths}) for each possible rank,
beginning with the rank of the Hamming code up to the full rank.
Since we did not use in the inductive step any propelinear codes
of lengths $2^{4s+3}-1,$ $s\geq 2$ and among them the propelinear
codes of lengths $63$, $127$ and $2^{11-1}$, this completes the
proof.

\medskip
{\bf Remarks}. In our opinion the open cases $n=r=63$ and
$n=r=127$ can be covered by the Vasil'ev construction applied to
full-rank propelinear perfect codes of lengths %% plural
 $31$ and $63$ using
a special propelinear functions. The last two open cases $n=127$,
$r=126$ and $n=r=2^{11-1}$ could be then
%%%# could we delete this "be"? No passive voice
 covered by the
Vasil'ev construction with the zero function $\lambda$ and by the
Mollard construction $\M(C^{2^4-1},C^{2^7-1})$ (or
$\M(C^{2^5-1},C^{2^6-1})$) with the zero function $f$
respectively.

The question of nontrivial lower and upper bounds on kernel
dimension, as well as the rank and kernel problem for propelinear
perfect codes are %% Џ®-¬®Ґ¬г, ¬­®¦ҐбвўҐ­­®Ґ зЁб«®:are as well as ўлбвгЇ Ґв Є Є б®о§ "and"
 still open. The rank and kernel problem can be
formulated as follows: which pairs of numbers $(r,k)$ are
attainable as the rank $r$ and kernel dimension $k$ of some
propelinear perfect code of length $n.$ Recall that the rank and
kernel problem for perfect binary codes was solved in \cite{AHS
rk}.

All computer searches have been carried out using the {\sc Magma}
\cite{MagmaCitation} software package. Some properties of perfect transitive codes of length $15$ and extended perfect transitive codes of length $16$ such as rank, dimension of the kernel, order of the automorphism group can be found in \cite{GusSolTablesArxiv}.

\medskip
{\bf Acknowledgement.} The authors  cordially thank Fedor Dudkin
for useful discussions.

\end{document}